\documentclass[10pt]{article}
\usepackage{repdefs,latexsym,longtable,graphicx,subeqn,varioref,cite,graphics}

\newcommand{\paperauthor}{Coralia Cartis\hspace*{-0.02cm}\footnotemark[1]\,, Nicholas I.~M.~Gould\footnotemark[2]\,\, and\,
 Philippe L.~Toint\footnotemark[3]}
\newcommand{\papertitle}{Universal regularization methods -- varying the power, the smoothness and the accuracy}

\title{\papertitle}
\author{\paperauthor}
\date{Revision completed August 19, 2018}

\newcommand{\theabstract}{}

\begin{document}
\renewcommand{\thefootnote}{\fnsymbol{footnote}}
\footnotetext[1]{Mathematical Institute, 
   Oxford University,
   Oxford OX2 6GG, UK.
   Email: coralia.cartis@maths.ox.ac.uk }
\footnotetext[2]{Computational Science and Engineering Department,
       Rutherford Appleton Laboratory, 
       Chilton, Oxfordshire, OX11 0QX, UK. 
       Email: nick.gould@stfc.ac.uk.}
\footnotetext[3]{
 NAXYS - University of Namur,
  61, rue de Bruxelles, B-5000, Namur, Belgium.
 Email: philippe.toint@unamur.be.}

\renewcommand{\thefootnote}{\arabic{footnote}}

\maketitle

\nocite{*}
\vspace*{-1cm}

\begin{abstract}
\theabstract
Adaptive cubic regularization methods have emerged as a credible alternative to linesearch and trust-region for 
smooth nonconvex optimization, with optimal complexity amongst second-order methods. Here we consider a general/new class of adaptive regularization methods,
that use first- or higher-order local Taylor models of the objective regularized by a(ny) power of the step size and applied to convexly-constrained optimization problems. We investigate the worst-case evaluation complexity/global rate of convergence
of these algorithms, when the level of sufficient smoothness  of the objective may be unknown or may even be absent. 
We find that the methods accurately 
reflect in their complexity  the degree of smoothness of
the objective and satisfy increasingly better bounds with improving accuracy of the models. The bounds vary continuously and
robustly with respect to the regularization power and accuracy of the model and the degree of smoothness of the objective. 
\end{abstract}

\vspace*{-0.45cm}
\setcounter{page}{1}
\numsection{Introduction}
\label{intro}

We consider the (possibly) convexly-constrained optimization problem
\eqn{UP}{
\min_{x\in \calF} f(x)
}
where $f:\Re^n\longrightarrow \Re$ is a smooth, possibly nonconvex, objective and where the feasible set $\calF\subset \Re^n$ is closed,
convex and non-empty (for example, the set $\calF$ could be described by simple bounds
and both polyhedral and more general convex constraints)\footnote{We are tacitly assuming that
the cost of evaluating constraint functions and their derivatives is negligible.}.  Clearly, the case of unconstrained optimization is covered here by letting
$\calF=\Re^n$. We are interested in the case when
$f\in \mathcal{C}^{p,\beta_p}(\calF)$, namely, $f$ is $p-$times continuously differentiable in $\calF$  with the $p$th derivative
being H\"older continuous of (unknown) degree $\beta_p \in [0,1]$\footnote{Note that if $\beta_p>1$, then the resulting class of objectives
is restricted to multivariate polynomials of degree $p$. If $p=1$, we only allow $\beta_1\in (0,1]$, for reasons to be explained later in the paper.}. 
We consider adaptive regularization methods applied to problem \req{UP} that 
generate feasible iterates $x_k$  that are (possibly very) approximate
minimizers over $\calF$ of local models of the form
\[
m_k(x_k+s)=T_p(x_k,s)+\frac{\sigma_k}{r}\|s\|_2^r,
\]
where $T_p(x_k,s)$ is the $p$th order Taylor polynomial of $f$ at $x_k$ and $r>p\geq 1$. 
The parameter $\sigma_k>0$ is adjusted to ensure
sufficient decrease in $f$ happens when the model value is decreased.
In this paper, we derive evaluation complexity bounds for finding
first-order critical points of \req{UP} using higher-order adaptive regularization methods. 
Despite the higher order of the models, the model minimization is performed only approximately, generalizing the approach in
\cite{bgmst}.  The proposed methods also ensure that the steps are `sufficiently long', in a new way, generalizing ideas in \cite{GN16}.
The ensuing complexity analysis shows the robust interplay of the regularization power $r$, the
model accuracy $p$ and the degree of smoothness  $\beta_p$ of the objective, with some surprising
results. In particular, we find that the degree of smoothness of the objective---which is often unknown and is even allowed to be absent here---is accurately reflected in the complexity
of the methods, independently of the regularization power, provided the latter is sufficiently large.
Furthermore, for all possible powers $r$, the methods satisfy increasingly better bounds
 as the accuracy $p$ of the models and  smoothness level $\beta_p$ are increased. All bounds vary continuously as a function of the regularization power and smoothness level.
Table 4.1 in Section 4 summarizes our complexity bounds.

We now review existing literature in detail and further clarify our approach, motivation 
and contributions. Cubic regularization for the (unconstrained) minimization of $f(x)$ for $x\in\Re^n$
was proposed independently by \cite{Grie81, NestPoly06, WeisDeufErdm07}, with \cite{NestPoly06} showing it has better global worst-case 
function evaluation complexity
than the method of steepest descent. Extending \cite{NestPoly06}, we proposed some practical variants -- Adaptive Regularization with Cubics (ARC) \cite{ARCpartII} -- that satisfy the same
complexity bound as the regularization methods in \cite{NestPoly06}, namely at most $\mathcal{O}(\epsilon^{-\frac{3}{2}})$ evaluations 
are needed to find a point $x$ for which  
\eqn{grad-eval}{
\|\nabla_x f(x)\|\leq \epsilon,
}  
under milder requirements on the algorithm (specifically, inexact model minimization). 
We further showed in \cite{cgt36, cgt38} that this complexity bound 
for ARC is sharp and optimal for a large class of second-order methods when applied to functions with globally Lipschitz-continuous second derivatives.  
Quadratic regularization, namely, a first order accurate model of the objective regularized by a quadratic term, has also been extensively studied, and shown to satisfy the complexity
bound of steepest descent, namely, $\mathcal{O}(\epsilon^{-2})$ evaluations to obtain \req{grad-eval} \cite{Nest04}. It was also shown in \cite{ARCpartII} that
one can loosen the requirement that global Lipschitz continuity of the second derivative holds, to just global H\"older continuity of the same derivative with
exponent $\beta_2\in (0,1]$.
Then, if one also
regularizes the quadratic objective model by the  power $2+\beta_2$ of the step, involving the (often unknown) H\"older exponent, the resulting method requires $\mathcal{O}(\epsilon^{-\frac{2+\beta_2}{1+\beta_2}})$ evaluations, which just as a function of $\epsilon$, belongs to the interval $\left[\epsilon^{-\frac{3}{2}}, \epsilon^{-2}\right]$; these bounds are sharp and optimal for objectives with corresponding level of smoothness of the Hessian \cite{cgt38}.  Note that this bound also holds if $\beta_2=0$.

An important related question and extension was answered in \cite{bgmst}:
if higher-order derivatives are available, can one improve the complexity of regularization methods? It was shown in \cite{bgmst}  that if one considers approximately minimizing
a $(r-1)$th order Taylor model of the objective regularized by the (weighted) $r$th power of the (Euclidean) norm of the step in each iteration (so $r=p+1$), the complexity of the resulting adaptive regularization 
method is $\mathcal{O}(\epsilon^{-\frac{r}{r-1}})$ evaluations to obtain \req{grad-eval}, under the assumption that the $(r-1)$th derivative tensor is globally Lipschitz continuous. The method proposed 
 in \cite{bgmst} measures progress of each iteration by comparing the Taylor model decrease (without the regularization term) to that of the true function decrease and only requiring mild approximate (local) minimization of the regularized model.  Here, we generalize these higher-order regularization methods from \cite{bgmst} to allow for an arbitrary local Taylor model, an arbitrary regularization power of the step and varying levels of 
smoothness of the highest-order derivative in the Taylor model. 


The interest in considering relaxations of Lipschitz continuity to H\"older continuity of derivatives comes not only from the needs of some engineering applications (such as flows in gas pipelines \cite[Section 17]{GPSA94} and properties of nonlinear
PDE problems \cite{BensFreh02}), but also in its own right in optimization theory, as a bridging case between the smooth and non-smooth classes of problems \cite{NemirovskiYudin, Nest13b}. In particular,  a zero H\"older exponent for a H\"older continuous derivative corresponds to a bounded derivative, 
an exponent in $(0,1)$ corresponds to a continuous but not necessarily  differentiable derivative, while an exponent of $1$ corresponds to a Lipschitz continuous derivative that can be differentiated again. For the case of function with H\"older-continuous gradients, methods have already been devised, and their complexity analysed,
both as a weaker set of assumptions and as an attempt to have a `smooth' transition between the smooth and nonsmooth (convex) problem classes, without knowing a priori the level of
smoothness of the gradient (i.e., the H\"older exponent) \cite{Devo13, Nest13b};
even lower complexity bounds are known \cite{NemirovskiYudin}. In \cite{cgt51} we considered regularization methods applied to nonconvex
objectives with H\"older continuous gradients (with unknown exponent $\beta_1\in (0,1]$), that employ a first-order quadratic model of the objective regularized by the $r$th power of the step.  
We showed that the worst-case complexity of the resulting regularization methods 
varies depending on $\min\{r,1+\beta_1\}$. In particular, when $1<r\leq 1+\beta_1$, the methods take at most
$\mathcal{O}\left(\epsilon^{-\frac{r}{r-1}}\right)$ evaluations/iterations until termination, and otherwise, at most $\mathcal{O}\left(\epsilon^{-\frac{1+\beta_1}{\beta_1}}\right)$ evaluations/iterations
to achieve the same condition. The latter complexity bound reflects the smoothness of the objective's landscape, without prior knowledge or use of it in the algorithm, and is independent of the
regularization power.
Here we generalize the approach in \cite{cgt51} to $p$th
order Taylor models and find that similar bounds can be obtained. Also, we are able to allow $\beta_p=0$ provided $p\geq 2$.
We note that advances beyond Lipschitz continuity of the derivatives for higher-order regularization methods were also obtained in
\cite{Chen}, 
where a class of problems with discontinuous and possibly infinite derivatives (such as when cusps are present) is analysed,
yielding similar  bounds to \cite{bgmst}.


Recently, \cite{GN16} proposed a new cubic regularization scheme that yields a {\it universal} algorithm in the sense that its complexity
 reflects the (possibly unknown or even absent) degree of sufficient smoothness of the objective; the approach in \cite{GN16}  addresses the case $p=2$, $r=3$ and $\beta_2\in [0,1]$ in our framework. 
Our ARp algorithm includes a modification in a similar (but not  identical) vein to that in \cite{GN16}.  In particular, our
approach
checks a theoretical condition that carefully monitors the length of the step
on each iteration on which the objective is sufficiently decreased. The technique in \cite{GN16} is different in that it requires a specific/new sufficient decrease condition of the objective on each iteration that makes progress. We generalize the approach in \cite{GN16} and achieve complexity bounds with similar universal properties  
for varying $r$, $p$ and unknown $\beta_p\in [0,1]$, provided $r\geq
p+\beta_p$. We are also able to analyze ARp's complexity in the regime $p<r\leq p+\beta_p$ providing continuously varying results with $r$ and $\beta_p$.

Our algorithm can be applied to convexly-constrained optimization problems with nonconvex objectives, where
the constraint/feasibility evaluations are inexpensive, offering another 
generalization of proposals in \cite{bgmst} and \cite{GN16} which are presented for the unconstrained case only; we also extend \cite{GN16}
by allowing inexact subproblem solution.

The structure of the paper is as follows. Section 2 describes our main algorithmic framework, ARp. Section 3 presents our complexity analysis while Section 4 concludes with a summary of our complexity bounds (see Table 4.1) and a discussion of the results.

\numsection{A universal adaptive regularization framework - ARp}

Let $f\in \mathcal{C}^p(\calF)$, with $p$ integer, $p\geq 1$; let $r\in \Re$, $r> p\geq 1$. 
We measure optimality using a suitable
continuous first-order criticality measure for \req{UP}.  We define this measure for a general
function $h:\Re^n\longrightarrow\Re$ on $\calF$: for an arbitrary
$x \in \calF$, the criticality measure is given by 
\beqn{omegadef}
\pi_h(x) \eqdef  \| P_\calF[x - \nabla_x h(x)] - x \|,
\eeqn
where $P_\calF$ denotes the orthogonal projection onto $\calF$ and $\|\cdot\|$
the Euclidean norm. Letting $h(x):=f(x)$ in \req{omegadef}, it is known that $x$ is a first-order critical point of
problem \req{UP} if and only if $\pi_f(x) = 0$. Also note that 
\[
\pi_f(x) =\|\nabla_x f(x)\| \tim {whenever $\calF = \Re^n$.}
\]
For more properties of this measure see
\cite{Bertsekas, TRbook}. 

Our ARp algorithm generates feasible iterates $x_k$ that (possibly very) approximately
minimize the local model
\eqn{2:localmodels}{
m_k(x_k+s)=T_p(x_k,s)+\frac{\sigma_k}{r}\|s\|^r \tim{subject to $x_k+s\in\calF$,}
}
which is a regularization of the $p$th order Taylor model of $f$ around $x_k$,
\eqn{2:qth-Taylor}{
T_p(x_k,s)=f(x_k)+\sum_{j=1}^p\frac{1}{j!}\nabla^j_x f(x^k)[s]^j,
}
where $\nabla^j_x f(x^k)[s]^j$ is the $j$th order tensor $\nabla^j_x f(x^k)$ of $f$ at $x_k$ applied to the vector $s$ repeated $j$ times. 
Note that $T_p(x_k,0)=f(x_k)$. We will also use the measure \req{omegadef} with $h(s):=m_k(x_k+s)$ for terminating the approximate minimization of 
$m_k(x_k+s)$, and for which we have again
\[
\pi_{m_k}(x_k+s) =\|\nabla_s m_k(x_k+s)\| \tim {whenever $\calF = \Re^n$.}
\]
A summary of the main algorithmic framework is as follows.
\algo{ARq-overest}{A universal ARp variant.}{
\vspace*{-0.3cm}
\begin{description}
\item[Step 0: Initialization.]
An initial point $x_0\in \calF$ and an initial regularization parameter
$\sigma_0\geq 0$ are given, as well as an accuracy level $\epsilon>0$.  The
constants  $\eta_1$,
$\eta_2$, $\gamma_1$, $\gamma_2$ and $\gamma_3$, $\theta$,
$\sigma_{\min}$
and $\alpha$, are also given and satisfy 
\eqn{unc1-arg-eta-gamma}{
 \theta>0,\quad 
\sigma_{\min}\in (0,\sigma_0],\quad 
0 < \eta_1 \leq \eta_2 < 1 \tim{and}
0<\gamma_3<1 < \gamma_1 < \gamma_2 \tim{and} \alpha\in \left(0,\frac{1}{3}\right].
}
Compute $f(x_0)$, $\nabla_x f(x_0)$ and set $k=0$. If $\pi_f (x_0)<\epsilon$,
terminate. Else, for $k\geq 0$, do:

\item[Step 1: Model set-up. ] 
Compute derivatives of $f$ of order $2$ to $p$ at $x_k$.

\item[Step 2: Step calculation. ]
Compute the step $s_k$ by approximately  minimizing the  model
$m_k(x_k+s)$ in \req{2:localmodels} over $x_k+s\in\calF$ such that the following conditions hold,
\eqn{feas-it}{x_k+s_k\in \calF,}
\eqn{2:str-decr}{
m_k(x_k+s_k) < f(x_k)
}
and 
\eqn{ho-cond}{
\pi_{m_k} (x_k+s_k)\leq \theta \|s_k\|^{r-1}.
}

\item[Step 3: Test for termination.] Compute $\nabla_x f(x_k+s_k)$. If $\pi_f (x_k+s_k) <\epsilon$, terminate
  with the approximate solution $x_\epsilon=x_k+s_k$.

\item[Step 4: Acceptance of the trial point. ]
Compute $f(x_k+s_k)$ and define
\eqn{rhokdef}{
\rho_k = \frac{f(x_k) - f(x_k+s_k)}{f(x_k) - T_p(x_k,s_k)}.
}
If $\rho_k \geq \eta_1$,  check whether
\eqn{step-long-gn}{
\sigma_k\|s_k\|^{r-1}\geq \alpha \pi_f (x_k+s_k).
}
If both $\rho_k \geq \eta_1$ and \req{step-long-gn} hold, then 
define
$x_{k+1} = x_k + s_k$; otherwise define $x_{k+1} = x_k$.

\item[Step 5: Regularization parameter update. ]
Set
\begin{equation}\label{sigupdate}
\sigma_{k+1} \in  \left\{\begin{array}{ll}
{}[\max(\sigma_{\min},\gamma_3\sigma_k), \sigma_k ]
& \tim{if} \rho_k \geq \eta_2\tim{and} \req{step-long-gn} \tim{holds,} \\
{}[\sigma_k, \gamma_1 \sigma_k ]                   & \tim{if} \rho_k
\in [\eta_1,\eta_2) \tim{and} \req{step-long-gn} \tim{ holds,}\\
{}[\gamma_1 \sigma_k, \gamma_2 \sigma_k ] & \tim{if} \rho_k < \eta_1
\tim{or} \req{step-long-gn} \tim{fails.}
  \end{array} \right.
  \end{equation}
Increment $k$ by one, and go to Step~1 
 if $\rho_k\geq \eta_1$ and \req{step-long-gn} hold, and to
Step 2 otherwise.
\end{description}
}

Iterations for which $\rho_k\geq \eta_1$ and \req{step-long-gn} hold (and so $x_{k+1}=x_k+s_k$) are called {\it successful}, 
those for which $\rho_k\geq \eta_2$ and \req{step-long-gn} hold are referred to as {\it very successful}, while the remaining ones are {\it unsuccessful}.
For  a(ny) $j\geq 0$, we denote the set of successful iterations up to $j$ 
by $\mathcal{S}_j=\{0\leq k\leq j:\,\rho_k\geq \eta_1 \tim{and} \req{step-long-gn} \,\,{\rm holds}\}$  and the set of unsuccessful ones by 
$\mathcal{U}_j=\{0,\ldots,j\}\setminus \mathcal{S}_j$.
We have the following simple lemma that relates the number of successful 
and unsuccessful iterations and that is ensured by the mechanism of the Algorithm \ref{ARq-overest}.


\llem{it-count-total-0}{\cite[Theorem 2.1]{ARCpartII} For any fixed $j\geq
  0$ until termination, let $\sigma_{\rm up}>0$ be such that $\sigma_k\leq \sigma_{\rm up}$
  for all $k\leq j$ in Algorithm \ref{ARq-overest}. Then 
\eqn{Uj}{
|\mathcal{U}_j|\leq \frac{|\log\gamma_3|}{\log\gamma_1}|\mathcal{S}_j|+\frac{1}{\log\gamma_1}\log \left(\frac{\sigma_{\rm up}}{\sigma_0}\right),
}
where $|\cdot|$ denotes the cardinality of the respective index set.}

\proof{The proof of \req{Uj} follows identically to the given reference; note that the sets $\mathcal{S}_j$ and $\mathcal{U}_j$ are not identical to the usual ARC
ones in \cite{ARCpartII} but the mechanism for modifying $\sigma_k$ in ARp coincides with the one in ARC on these iterations and that is why the proof
of this lemma follows identically to \cite[Theorem 2.1]{ARCpartII}.}


Now we comment on the construction of the ARp algorithm. Note that
the model minimization conditions (Step 2) and the definition of $\rho$ in Step 4
are straightforward generalizations of the approach in \cite{bgmst} to $p$th order Taylor models regularized by
different powers  $r$ of the norm of the step. Furthermore, recall that conditions \req{feas-it},
\req{2:str-decr} and \req{ho-cond} are approximate {\it local} optimality conditions for the nonconvex polynomial model $m_k(x_k+s)$
minimization over a convex set, $x_k+s\in \mathcal{F}$; in fact, they are even weaker than that as they require 
 strict decrease (from the base point $s=0$) and {\it approximate first-order criticality} for the convexly constrained model.
Thus, any descent optimization method---even first-order algorithms such as the projected
gradient method---can be applied to ensure these conditions with ease (with no additional derivatives evaluations required than those needed to set up the model $m_k$ at $x_k$). Designing efficient techniques specifically for the approximate minimization of such regularized, nonconvex, high-order polynomial optimization problems is beyond our scope here, but an essential component of the success of such methods.  Existing regularization-related approaches are available for general nonconvex problems up to third order \cite{bgms1, bgms2}, or 
dedicated to convex regularized tensor models (see \cite{NestModels} and the references therein) or specialized to nonlinear
least-squares problems \cite{Gould1, Gould2}; these complement classical references such as \cite{Schnabel}, where
 third and fourth
order tensor methods were proposed.

However, there are two main differences to the by-now standard 
approaches to (cubic or higher order) regularization methods. Firstly, we check whether the gradient goes below $\epsilon$ at each trial points, and if so, terminate
on possibly unsuccessful iterations (Step 3). Secondly, when the step $s_k$ provides sufficient decrease according to \req{rhokdef},
we check whether $s_k$ satisfies \req{step-long-gn}, and only allow  steps that have such carefully-monitored length to be taken by the algorithm; if \req{step-long-gn} fails  or $\rho_k\leq \eta_1$, $\sigma_k$ is increased. 
Note that though the length of the step $s_k$ decreases as $\sigma_k$ is increased, this is not the case for the expression $\sigma_k\|s_k\|^{r-1}$ in \req{step-long-gn}, which increases with $\sigma_k$, as  Lemma \ref{lem-sigma-k-general} implies. 
These two additional ingredients---the gradient calculation at each trial point and the step length condition \req{step-long-gn}---are directly related to trying to achieve universality of ARp, extending ideas from \cite{GN16}. 
Further explanations and discussions for the theoretical need, or otherwise, for condition \req{step-long-gn} are given 
next, in Remark \ref{remark1}, and later in the paper, in Remarks \ref{remark2} (b) and  \ref{remark4} (b).  

\begin{remark}\label{remark1}
{\rm We further comment on condition \req{step-long-gn}, its connections to \cite{GN16} and existing literature, and possible alternatives.
\begin{itemize}
\item[(a)] We can replace condition \req{step-long-gn} with the weaker requirement that
$\sigma_k\|s_k\|^{r-1}\geq \alpha \epsilon$; then, all subsequent results would remain unchanged. This choice however, would make 
the algorithm construction dependent on the accuracy $\epsilon$ (elsewhere than in the termination condition), 
which is not numerically advisable. 

\item[(b)] Instead of requiring \req{step-long-gn} on each successful step, we could ask that each model minimization step calculated
in Step 2 satisfies \req{step-long-gn}; 
if \req{step-long-gn} failed, $\sigma_k$ would be increased at the end of Step 2 and the model minimization step would be repeated. 
This approach may result in an unnecessarily small step in practice, but the ensuing ARp complexity bounds would remain qualitatively similar.

\item[(c)] 
Condition \req{step-long-gn} does not appear as such in the algorithmic variants proposed in \cite{GN16}, as those enforce
sufficient decrease conditions on $f$ in the algorithm for the case $p=2$ and $r=3$, which is the only case addressed in \cite{GN16}. But \req{step-long-gn}
(with $r=3$) is a necessary ingredient for achieving the required sufficient decrease conditions
in \cite{GN16}; see Lemma 2.3 (in particular, equation (2.21)) therein.





\item[(d)]  Following \cite{GN16}, instead of \req{step-long-gn}, we could employ a different definition of 
$\rho_k$ in \req{rhokdef},
 namely, replacing the denominator in \req{rhokdef}  by a rational function in $\epsilon$ and $\sigma_k$, or by 
a function of $\sigma_k$ and the gradient at the new point (see for example \cite[(6.5)]{GN16}),
to achieve the desired order of model/function decrease for universal complexity and behaviour. 
According to our calculations, again, qualitatively similar complexity bounds would be obtained for such ARp variants. 

We note that using specific $\rho_k$ definitions (namely, with a denominator connected to the length of the step)
 so as to enforce a particular sufficient decrease property for the objective evaluations
was also used in \cite{TRACE, Birgin} for trust-region and quadratic regularization variants, in order  to achieve optimal complexity bounds
for the ensuing methods.

\item[(e)] According to our calculations, without the condition 
\req{step-long-gn} on the length of the step, or a similar measure of progress, the complexity of ARp would   dramatically (but continuously) worsen 
 in the regime when $r>p+\beta_p$, as $r$ increases. But as we clarify at the end of  Section \ref{complexity-ARP-section}, for the case $r\leq p+\beta_p$, same-order complexity bounds could be obtained for ARp without using \req{step-long-gn}; so in
  principle, for this parameter regime, \req{step-long-gn} could be removed from the construction of ARp. 
 However, note that as $\beta_p$ is not generally known {\it a priori}, the regime of most interest 
-- both in terms of best complexity bounds and practicality --  is when 
$r$ is large; hence the need for condition \req{step-long-gn} in ARp, for both regimes.

\end{itemize}}
\end{remark}

\numsection{Worst-case complexity analysis of ARp}
\label{complexity-ARP-section}

\subsection{Some preliminary properties}

We have the following simple consequence of \req{2:str-decr}.

\llem{Lemma-decrease}{On each iteration of Algorithm
  \ref{ARq-overest},  we have the decrease
\eqn{2:model-decr}{
f(x_k)-T_p(x_k,s_k)\geq \frac{\sigma_k}{r}\|s_k\|^r.
}
}

\proof{Note that condition \req{2:str-decr} and the definition of
  $m_k(s)$ in \req{2:localmodels} immediately give
  \req{2:model-decr}.}

We have the following upper bound on $s_k$.

\llem{lem-upper-bound-sk}{
On each iteration of Algorithm
  \ref{ARq-overest}, we have
\eqn{sk-upper-bd}{
\|s_k\|\leq 
\max_{1\leq j\leq p}\left\{\left(\frac{pr}{j!\sigma_k}\|\nabla^j_x f(x_k)\|\right)^{\frac{1}{r-j}}\right\}.
}}

\proof{It follows from \req{2:str-decr}, \req{2:localmodels} and
  \req{2:qth-Taylor} that 
\[s_k^T\nabla_x f(x_k)+\frac{1}{2}\nabla^2_x
f(x_k)[s_k,s_k]+\ldots+\frac{1}{p!}\nabla^p_x
f(x^k)[s_k,s_k,\ldots,s_k]+\frac{\sigma_k}{r}\|s_k\|^r<0,
\]
which from Cauchy-Schwarz and norm properties, further implies
\[
-\|s_k\|\cdot\|\nabla_x f(x_k)\|-\frac{1}{2}\|s_k\|^2\cdot\|\nabla^2_x
f(x_k)\|-\ldots-\frac{1}{p!}\|s_k\|^p\cdot\|\nabla^p_x
f(x^k)\|+\frac{\sigma_k}{r}\|s_k\|^r<0,
\]
or equivalently,
\[
\sum_{j=1}^p \left(\frac{\sigma_k}{pr}\|s_k\|^r-\frac{1}{j!}\|s_k\|^j\cdot\|\nabla^j_x
f(x^k)\|\right)<0.
\]
The last displayed equation cannot hold unless at least one of the
terms on the left-hand side is negative, which is equivalent to
\req{sk-upper-bd}, using also that $r>p\geq 1$.}


Let us assume that $f\in\mathcal{C}^{p,\beta_p}$, namely,
\begin{description}
\item[A.1] $f\in C^p(\calF)$ and $\nabla^p_x f$ is 
H\"older continuous on the path of the iterates and trial points, namely,
\[
\|\nabla^p_x f(y)-\nabla^p_x f(x_k)\|_T\leq (p-1)! L_p\|y-x_k\|^{\beta_p}
\]
 holds for all $y \in [x_k,x_k+s_k]$, $k\geq 0$ and some constants $L_p \geq 0$ and
 $\beta_p \in [0,1]$, where $\|\cdot\|$ is the Euclidean norm on
 $\Re^n$ and $\|\cdot\|_T$ is recursively induced by this norm on the
 space of the $p$th order tensors. 
\end{description}

A simple consequence of A.1 is that
\eqn{approx}{
|f(x_k+s_k) - T_p(x_k,s_k)|\leq
\frac{L_p}{p}\|s_k\|^{p+\beta_p},\quad k\geq 0,
}
and
\eqn{approx1}{
\|\nabla_x f(x_k+s_k) - \nabla_s T_p(x_k,s_k)\|\leq
L_p\|s_k\|^{p+\beta_p-1},\quad k\geq 0;
}
see \cite{bgmst} for a proof of \req{approx} and \req{approx1}, with A.1 replacing Lipschitz continuity of the $p$th derivative.

\begin{remark}\label{remark5} 
{\rm Note that throughout the paper we assume $r>p\geq 1$, $r\in\Re$ and $p\in \Na$; and that either $p\geq 1$ and $\beta_p\in (0,1]$
or $p\geq 2$ and $\beta_p\in [0,1]$. Thus in both cases $p+\beta_p-1>0$.}
\end{remark}

Two useful preliminary lemmas follow.

\llem{lem-step-lowerbd}{Assume that A.1 holds. Then on each iteration
  of Algorithm \ref{ARq-overest}, we have
\eqn{step-lowerbd}{
\pi_f(x_k+s_k)\leq L_p\|s_k\|^{p+\beta_p-1}+(\sigma_k+\theta)\|s_k\|^{r-1}.
}}

\proof{Using the triangle inequality and \req{omegadef} with $h\eqdef f$ and $h\eqdef m_k$, we obtain 
\[
\begin{array}{lcl}
\pi_f(x_k+s_k)
& = & \| 
P_\calF[x_k+s_k-\nabla_x f(x_k+s_k) ] - P_\calF[x_k+s_k-\nabla_s m_k(x_{k}+s_k) ]
\\*[2ex]
& & \;\; + \; P_\calF[x_k+s_k-\nabla_s m_k(x_k+s_k) ] - (x_k+s_k) \|
\\*[2ex]
& \leq &  \|P_\calF[x_k+s_k-\nabla_x f(x_k+s_k) ] -  P_\calF[x_k+s_k-\nabla_s m_k(x_{k}+s_k) ] \| 
              + \pi_{m_k}(x_{k}+s_k). \\*[2ex]
\end{array}
\]
The last inequality, the contractive property of the projection operator $P_{\calF}$ and 
 the  inner termination condition \req{ho-cond} give
\eqn{cgt54:temp-3}{ 
\pi_f(x_k+s_k)\leq   \| \nabla_x f(x_k+s_k)  -  \nabla_s m_k(x_k+s_k) \|  + \theta \|s_k\|^{r-1}.
}
We have from \req{2:localmodels} that 
\[
\nabla_s m_k(x_k+s) =\nabla_s T_p(x_k,s) + \sigma_k \|s\|^{r-1}\frac{s}{\|s\|}
\]
and so 
\eqn{temp-bd}{
\begin{array}{lcl}
\|\nabla_x f(x_k+s_k)-\nabla_s m_k(x_k+s_k)\| &\leq& \|\nabla_x f(x_k+s_k)-\nabla_s T_p(x_k,s_k)\| +
\sigma_k\|s_k\|^{r-1} \\*[2ex]
&\leq& L_p\|s_k\|^{p+\beta_p-1} +\sigma_k\|s_k\|^{r-1} ,
\end{array}
}
where we used \req{approx1} to obtain the second inequality.  
Now \req{step-lowerbd} follows from replacing \req{temp-bd} in \req{cgt54:temp-3}.}

\llem{lem-sigma-k-general}{Assume that A.1 holds. If
\eqn{sigma-general}{\sigma_k\geq \max\left\{\theta, 
  \kappa_2\|s_k\|^{p+\beta_p-r} \right\},}
where 
\eqn{constant-0}{\kappa_2\eqdef \frac{rL_p}{p(1-\eta_2)},}
then both $\rho_k\geq \eta_2$ and \req{step-long-gn} hold, and so iteration $k$ is
very successful.}

\proof{We assume that \req{sigma-general} holds, which implies that
\eqn{cgt54:temp-sigma}{\sigma_k\geq \kappa_2\|s_k\|^{p+\beta_p-r}.}
The definition of $\rho_k$ in \req{rhokdef} gives
$\displaystyle |\rho_k-1|=
\frac{|f(x_k+s_k)-T_p(x_k,s_k)|}{f(x_k)-T_p(x_k,s_k)}$,
whose numerator we upper bound by \req{approx}, and whose denominator we lower bound by \req{2:model-decr}, 
to deduce
\eqn{cgt54:temp}{
|\rho_k-1|
\leq\frac{\frac{L_p}{p} \|s_k\|^{p+\beta_p}}{\frac{\sigma_k}{r}\|s_k\|^r}=\frac{rL_p}{p\sigma_k}\|s_k\|^{p+\beta_p-r}.
}
We employ \req{cgt54:temp-sigma} and the expression of $\kappa_2$ 
in \req{constant-0},  in \req{cgt54:temp}, to deduce that 
 $|1-\rho_k|\leq 1-\eta_2$, which ensures that $\rho_k\geq \eta_2$.

It remains to show that  \req{sigma-general} also implies \req{step-long-gn}.
From \req{sigma-general}, we have that $\sigma_k\geq \theta$, which together with
 \req{step-lowerbd}, give
\eqn{cgt54:temp-1}{\pi_f(x_k+s_k)\leq \|s_k\|^{p+\beta_p-1} \left(L_p+2\sigma_k\|s_k\|^{r-p-\beta_p}\right).}
The definition \req{constant-0}, 
and requirements $r>p$ and $\eta_2\in (0,1)$, imply that
 $L_p\leq \kappa_2$. This and \req{cgt54:temp-1} give
\eqn{cgt54:temp-2}{\pi_f(x_k+s_k)\leq \|s_k\|^{p+\beta_p-1} \left(\kappa_2+2\sigma_k\|s_k\|^{r-p-\beta_p}\right).}
From \req{cgt54:temp-sigma}, $\kappa_2\leq \sigma_k \|s_k\|^{r-p-\beta_p}$. We use this to bound 
 $\kappa_2$ in \req{cgt54:temp-2}, which gives the inequality 
\[
\pi_f(x_k+s_k)\leq \|s_k\|^{p+\beta_p-1} \left(3\sigma_k\|s_k\|^{r-p-\beta_p}\right)=3\sigma_k \|s_k\|^{r-1}.
\]
Thus $\sigma_k\|s_k\|^{r-1}\geq \frac{1}{3}\pi_f(x_k+s_k)$, which implies \req{step-long-gn} since $\alpha\leq \frac{1}{3}$.} 


\subsection{The case when $r> p+\beta_p$}

Using Lemmas \ref{lem-step-lowerbd} and \ref{lem-sigma-k-general}, we
have the following result, which together with its proof,  were inspired by and generalize the result and proof in  \cite[Lemma 2.3]{GN16}.

\llem{lem-gn-sigma-success}{Let $r> p+\beta_p$  and assume A.1.  While Algorithm
  \ref{ARq-overest} has not terminated, if
\eqn{gn-sigma-lower}{\sigma_k\geq \max\left\{\theta, \kappa_1 \epsilon^{\frac{p+\beta_p-r}{p+\beta_p-1}} \right\},}
where 
\eqn{constants}{\kappa_1\eqdef \left(3^{r-p-\beta_p}\kappa_2^{r-1}\right)^{\frac{1}{p+\beta_p-1}} \tim{and} \kappa_2 \tim{is defined in \req{constant-0},}}
then  \req{sigma-general}  holds, and so iteration $k$ is very successful.}

\proof{
We will prove our result by contradiction. We assume that \req{sigma-general} does not hold on iteration $k$, and so 
\eqn{sigma-g-1}{\sigma_k\|s_k\|^{r-p-\beta_p}<\kappa_2.}
Note that while Algorithm \ref{ARq-overest} does not terminate, we have $\pi_f(x_k+s_k)\geq \epsilon$. Also, from 
\req{gn-sigma-lower}, $\sigma_k\geq \theta$. We use these two inequalities into \req{step-lowerbd} to deduce
\eqn{cgt54:temp-4}{
\epsilon\leq L_p\|s_k\|^{p+\beta_p-1} +2\sigma_k\|s_k\|^{r-1}= \|s_k\|^{p+\beta_p-1} \left(L_p+2\sigma_k\|s_k\|^{r-p-\beta_p}\right).
}
We now employ \req{sigma-g-1} to upper bound the second term in \req{cgt54:temp-4} by $2\kappa_2$, namely,
\eqn{cgt54:temp-5}{
\epsilon< \|s_k\|^{p+\beta_p-1} \left(L_p+2\kappa_2\right).
}
We use \req{sigma-g-1} again to provide an upper bound on $\|s_k\|$, which is possible since $r>p+\beta_p$. Thus
\eqn{s:bd-r}{
\|s_k\|\leq \left(\frac{\kappa_2}{\sigma_k}\right)^{\frac{1}{r-p-\beta_p}}.
}
Using this bound in \req{cgt54:temp-5}, which is possible since $p+\beta_p>1$, we obtain the first inequality below,
\eqn{cgt54:temp-6}{
\epsilon<\left(\frac{\kappa_2}{\sigma_k}\right)^{\frac{p+\beta_p-1}{r-p-\beta_p}}\left(L_p+2\kappa_2\right)
< \left(\frac{\kappa_2}{\sigma_k}\right)^{\frac{p+\beta_p-1}{r-p-\beta_p}} \cdot (3\kappa_2),
}
where to obtain the second inequality, we used that $L_p< \kappa_2$, which in turn follows from
\req{constant-0}, $r>p$ and $\eta_2\in (0,1)$.
Finally, \req{cgt54:temp-6}  and the definition of $\kappa_1$ in \req{constants}
 imply that  $\sigma_k<\kappa_1 \epsilon^{\frac{p+\beta_p-r}{p+\beta_p-1}}$, which contradicts \req{gn-sigma-lower}.
Thus \req{sigma-general} must hold and Lemma \ref{lem-sigma-k-general} implies that $\rho_k\geq \eta_2$ and \req{step-long-gn} hold, and so $k$ is very successful.} 

\begin{remark}\label{remark2}
{\rm \begin{itemize}
\item[(a)] (Parameter regime)
The proof of Lemma \ref{lem-gn-sigma-success} requires $r>p+\beta_p$ and $p+\beta_p>1$ 
(to deduce \req{s:bd-r} and \req{cgt54:temp-6}, respectively). However, the result of Lemma \ref{lem-gn-sigma-success}
remains true if $r=p+\beta_p$ and it is proved together with the case $r<p+\beta_p$ in Lemma \ref{p-small:sigma}.  Note that, when 
$r=p+\beta_p$, \req{gn-sigma-lower} becomes $\sigma_k\geq \max\{\theta, \kappa_2\}$, which precisely
matches the corresponding expression \req{good-bd} in Lemma \ref{p-small:sigma} for this same case.

\item[(b)] (Condition \req{step-long-gn}) 
Without employing \req{step-long-gn}, 
we showed  inequality \req{step-lowerbd} that connects the length of the step to that of the projected gradient. The two terms on the right-hand side of \req{step-lowerbd} have similar forms as powers of $\|s_k\|$, with the exponents crucially determined by H\"older continuity properties of the objective and the power of the regularization term in the model, respectively. 
Lemmas \ref{lem-sigma-k-general} and \ref{lem-gn-sigma-success}  proved that if $\sigma_k$ is sufficiently large, then the second term in 
\req{step-lowerbd}, namely, $\sigma_k\|s_k\|^{r-1}$, will be larger than the term that is a multiple of  $\|s_k\|^{p+\beta_p-1}$;
hence ensuring that \req{step-long-gn} holds. 
To further explain this point, note that 
in \req{step-lowerbd}, when  $r>p+\beta_p$ and $\|s_k\|\leq 1$ (which is the difficult case), the larger term on the 
right-hand side is a multiple of
$\|s_k\|^{p+\beta_p-1}$ when $\sigma_k$ is larger than a constant.
  Lemma \ref{lem-gn-sigma-success} showed that if $\sigma_k$ is further increased,
in an $\epsilon$-dependent way, then the term that is a multiple of $\|s_k\|^{r-1}$
in  \req{step-lowerbd} becomes the larger of the two terms.

\end{itemize}}
\end{remark}

\llem{lem-sigma-upper-bd}{Let $r> p+\beta_p$  and assume A.1.  Then, while Algorithm
  \ref{ARq-overest} has not terminated, we have
\eqn{sigma-upper}{\sigma_k\leq \max\left\{\sigma_0,\gamma_2\theta, \gamma_2\kappa_1 \epsilon^{\frac{p+\beta_p-r}{p+\beta_p-1}} \right\},}
where $\kappa_1$ is defined in \req{constants}. }

\proof{Let the right-hand side of \req{gn-sigma-lower} be denoted by $\overline{\sigma}$.
It follows from Lemma \ref{lem-gn-sigma-success} and the mechanism of the algorithm that
\eqn{implication}{
 \sigma_k\geq \overline{\sigma} \quad\Longrightarrow\quad \sigma_{k+1}\leq \sigma_k. 
}
Thus, when $\sigma_0\leq \gamma_2\overline{\sigma}$,  it follows that $\sigma_k\leq \gamma_2\overline{\sigma}$,
where the factor $\gamma_2$ is introduced for the case when $\sigma_k$ is less than $\overline{\sigma}$ and the iteration $k$ is not very successful.  
Letting $k=0$ in \req{implication} gives \req{sigma-upper} when $\sigma_0\geq \gamma_2\overline{\sigma}$ since $\gamma_2>1$.
}

We are ready to establish an upper bound on the number of successful iterations until termination.

\lthm{thm-success}{Let $r> p+\beta_p$, assume A.1 and that $\{f(x_k)\}$ is bounded below by $f_{\rm low}$ and $\epsilon\in (0,1]$. Then for all successful iterations $k$ until the termination of
Algorithm  \ref{ARq-overest}, we have
\eqn{success-k}{f(x_k)-f(x_{k+1})\geq \kappa_{s,p} \epsilon^{\frac{p+\beta_p}{p+\beta_p-1}},
}
where 
\eqn{constants-success}{
\kappa_{s,p}\eqdef \frac{\eta_1}{r}\left(\frac{\alpha^{r}}{\sigma_{\max}}\right)^{\frac{1}{r-1}},\quad \sigma_{\max}\eqdef \max\left\{\sigma_0,\gamma_2\theta, \gamma_2\kappa_1\right\},
}
and $\kappa_1$ is defined in \req{constants}.
Thus Algorithm \ref{ARq-overest} takes at most 
\eqn{succ-no}{
\left\lfloor \frac{f(x_0) - f_{\rm low}}{\kappa_{s,p}} \epsilon^{-\frac{p+\beta_p}{p+\beta_p-1}}\right\rfloor
}
successful iterations/evaluations of derivatives of degree $2$ and above of $f$ until termination.
}

\proof{On every successful iteration $k$, we have $\rho_k\geq \eta_1$; this and Lemma \ref{Lemma-decrease} imply
\eqn{plarge:fct-d}{
f(x_k)-f(x_{k+1})\geq \eta_1(f(x_k)-T_p(x_k,s_k))\geq \eta_1 \frac{\sigma_k}{r}\|s_k\|^r=\frac{\eta_1}{r} (\sigma_k\|s_k\|^{r-1})\|s_k\|.}
On every successful iteration $k$ we also have that \req{step-long-gn} holds. Thus, while the algorithm has not terminated, we have
\eqn{cgt54:temp-7}{\sigma_k\|s_k\|^{r-1}\geq \alpha\epsilon \quad{\rm and}\quad \|s_k\|\geq \left(\frac{\alpha\epsilon}{\sigma_k}\right)^{\frac{1}{r-1}}.}
Applying the first and then the second inequality in \req{cgt54:temp-7} into \req{plarge:fct-d}, we deduce
\eqn{plarge:fct-d2}{
f(x_k)-f(x_{k+1})\geq \frac{\eta_1}{r}\alpha\epsilon\|s_k\|\geq \frac{\eta_1}{r}\alpha\epsilon\left(\frac{\alpha\epsilon}{\sigma_k}\right)^{\frac{1}{r-1}}=
\frac{\eta_1}{r}\frac{(\alpha\epsilon)^{\frac{r}{r-1}}}{\sigma_k^{\frac{1}{r-1}}}.}
We use that $\epsilon\in (0,1]$ in \req{sigma-upper} to deduce that
\eqn{sigma-useful}{
\sigma_k\leq \sigma_{\max}\epsilon^{\frac{p+\beta_p-r}{p+\beta_p-1}},
}
where $\sigma_{\max}$ is defined in \req{constants-success}. We combine this upper bound with \req{plarge:fct-d2} to see that
\[
f(x_k)-f(x_{k+1})\geq \frac{\eta_1}{r}(\alpha\epsilon)^{\frac{r}{r-1}}\sigma_{\max}^{-\frac{1}{r-1}} \epsilon^{\frac{r-p-\beta_p}{(p+\beta_p-1)(r-1)}}=\frac{\eta_1}{r}\left(\frac{\alpha^{r}}{\sigma_{\max}}\right)^{\frac{1}{r-1}}
\cdot \epsilon^{\frac{p+\beta_p}{p+\beta_p-1}},
\] 
which gives \req{success-k}. Using that $f(x_k)=f(x_{k+1})$ on unsuccessful iterations, and that $f(x_k)\geq f_{\rm low}$ for all $k$, we can sum up over all successful iterations to deduce 
\req{succ-no}.}



We are left with counting the number of unsuccessful iterations until termination, and the total iteration and evaluation upper bound.

\llem{it-count-total}{Let $r> p+\beta_p$ and $\epsilon\in (0,1]$. Then, for any fixed $j\geq
  0$ until termination, Algorithm \ref{ARq-overest} satisfies
\eqn{Uj-2}{
|\mathcal{U}_j|\leq \frac{|\log\gamma_3|}{\log\gamma_1}|\mathcal{S}_j|+\frac{1}{\log\gamma_1}\log \frac{\sigma_{\rm max}}{\sigma_0}
+ \frac{r-p-\beta_p}{(p+\beta_p-1)\log\gamma_1}|\log \epsilon|,}
where $\sigma_{\max}$ is defined in \req{constants-success}.}

\proof{We apply Lemma \ref{it-count-total-0}. 
To prove \req{Uj-2}, we use $\epsilon \in (0,1]$ and the upper bound \req{sigma-useful} in place of $\sigma_{\rm up}$ in \req{Uj}.}

\lcor{cor-complexity}{Let $r> p+\beta_p$  and assume A.1, that $\{f(x_k)\}$ is bounded below by $f_{\rm low}$ and  $\epsilon\in (0,1]$. 
Then Algorithm \ref{ARq-overest} takes at most 
\eqn{total-no}{
\left\lfloor \frac{f(x_0) - f_{\rm low}}{\kappa_{s,p}}\left(1+ \frac{|\log\gamma_3|}{\log\gamma_1}\right) \epsilon^{-\frac{p+\beta_p}{p+\beta_p-1}} 
+ \frac{r-p-\beta_p}{(p+\beta_p-1)\log\gamma_1}|\log \epsilon| + \frac{1}{\log\gamma_1}\log \frac{\sigma_{\rm max}}{\sigma_0}\right\rfloor
}
iterations/evaluations of $f$ and its derivatives until termination, where $\kappa_{s,p}$ and $\sigma_{\max}$ are defined in \req{constants-success}.}

\proof{The proof follows from Theorem \ref{thm-success} and \req{Uj-2}, where we let $j$ denote the first iteration with $\pi_f(x_j+s_j)<\epsilon$
(so the iteration where ARp terminates) and we use $j=|\mathcal{S}_j|+|\mathcal{U}_j|$.}


\begin{remark}\label{remark3}
{\rm \begin{itemize}
\item[(a)] (Comment on $\sigma_{\min}$)
We  note that the lower bound on $\sigma_k$,  $\sigma_k\geq \sigma_{\min}\geq 0$ for all $k$, imposed in \req{sigupdate}, has not been employed in 
the above proofs and it is also not needed when $r=p+\beta_p$. It seems that in the case
$r\geq p+\beta_p$, such a lower bound on $\sigma_k$ may follow implicitly from \req{step-long-gn}.
However, the requirement involving $\sigma_{\min}>0$ is needed for the case $r< p+\beta_p$.

\item[(b)] (Comment on $\epsilon$) In our main complexity results (such as Corollary \ref{cor-complexity}), we have a restriction
on the required accuracy tolerance $\epsilon \in (0,1]$; this restriction is for simplicity and simplification of expressions, so as to
capture dominating terms in the complexity bounds.  It is also intuitive, as we think of $\epsilon$ as (arbitrarily)
`small' compared to problem 
constants. Indeed, instead of an upper bound of $1$ on $\epsilon$, we could have used a bound depending on problem constants
such as $L_p$,  which would preserve the same dominating terms in the complexity bounds. However, as most such problem constants are generally unknown, we prefer our approach as it gives the users/readers a concrete value they can use. 

\end{itemize}}
\end{remark}

The constants in the bound \req{total-no} and their behaviour with respect to increasing values of $p$ 
are discussed in Section \ref{constants-page}.



\subsection{The case when $p<r\leq p+\beta_p$}

Note that $p<r\leq p+\beta_p$ imposes that $\beta_p>0$ in this case. Also, note that  
the proof of Lemma \ref{lem-gn-sigma-success} fails to hold for $r\leq p+\beta_p$. Thus we need a different approach here
to upper bounding $\sigma_k$. In particular, we need the following additional assumption (for the case when $r<p+\beta_p$).

\begin{description}
\item[A.2]  For $j\in\{1,\ldots,p\}$, the derivative $\{\nabla^j f(x_k)\}$
  is uniformly bounded above with respect to $k$, namely,
\[
\|\nabla^j f(x_k)\|\leq M_j \tim{for all $k\geq 0$,\quad
  $j\in\{1,\ldots,p\}$. }
\]
We let $M\eqdef \displaystyle\max_{1\leq j\leq p}\left\{\left(\frac{rp}{j!\sigma_{\min}}M_j\right)^{\frac{1}{r-j}}\right\}$
where $\sigma_{\min}$ is defined in \req{sigupdate}.
\end{description}

\llem{p-small:sigma}{Let $r\leq  p+\beta_p$  and assume A.1.  If $r<p+\beta_p$ assume also A.2 and $\sigma_{\min}>0$. 
  If
\eqn{good-bd}{\sigma_k\geq \max\left\{\theta, \kappa_2M^{p+\beta_p-r}\right\},}
where $\kappa_2$ and $M$ are defined in \req{constant-0} and A.2, respectively, then
\req{sigma-general}  holds, and so iteration $k$ is very successful.}

\proof{If $r=p+\beta_p$, then \req{good-bd} clearly implies \req{sigma-general} and so Lemma \ref{lem-sigma-k-general} applies.

If $r<p+\beta_p$, then we upper bound 
$\|s_k\|$ by using A.2 in \req{sk-upper-bd}, as well as $\sigma_k\geq
\sigma_{\min}$, to deduce that $\|s_k\|\leq M$ where $M$ is defined in A.2.
Now \req{good-bd} implies \req{sigma-general} and so Lemma \ref{lem-sigma-k-general} again applies, yielding that iteration $k$ is very successful.}

We are ready to  bound $\sigma_k$ from above for all iterations.

\llem{p:small-s-upper-bd}{Let $r\leq p+\beta_p$  and assume A.1. 
 If $r<p+\beta_p$ assume also A.2 and $\sigma_{\min}>0$. While Algorithm
  \ref{ARq-overest} has not terminated, we have
\eqn{psmall:sigma-upper}{\sigma_k\leq \max\left\{\sigma_0,\gamma_2\theta, \gamma_2 \kappa_2M^{p+\beta_p-r} \right\}\eqdef \sigma_{\rm up},}
where   $\kappa_2$ and $M$ are defined in \req{constant-0} and A.2, respectively.}

\proof{The proof follows a similar argument to that of Lemma \ref{lem-sigma-upper-bd}, with \req{gn-sigma-lower} replaced by \req{good-bd}. Note also that as $\epsilon$ does not appear in the bound \req{good-bd}, \req{psmall:sigma-upper} yields a constant
upper bound on $\sigma_k$ that is valid for all $k$, irrespective of the required accuracy level $\epsilon$.}
We are now ready to upper bound the number of successful iterations of Algorithm \ref{ARq-overest} until termination.

\lthm{p-small:thm-success}{Let $r\leq p+\beta_p$, assume A.1 and that $\{f(x_k)\}$ is bounded below by $f_{\rm low}$. 
 If $r<p+\beta_p$ assume also A.2 and $\sigma_{\min}>0$.
Then for all successful iterations $k$ until the termination of
Algorithm  \ref{ARq-overest}, we have
\eqn{psmall:fct-final}{f(x_k)-f(x_{k+1})\geq \kappa_{s,r} \epsilon^{\frac{r}{r-1}},}
where 
\eqn{psmall:constants-success}{
\kappa_{s,r}\eqdef \frac{\eta_1}{r}\left(\frac{\alpha^{r}}{\sigma_{\rm up}}\right)^{\frac{1}{r-1}},}
and $\sigma_{\rm up}$ is defined in \req{psmall:sigma-upper}.
Thus Algorithm \ref{ARq-overest} takes at most 
\eqn{psmall:succ-no}{
\left\lfloor \frac{f(x_0) - f_{\rm low}}{\kappa_{s,r}} \epsilon^{-\frac{r}{r-1}}\right\rfloor
}
successful iterations/evaluations of derivatives of degree $2$ and higher of $f$ until termination.}

\proof{Note that \req{plarge:fct-d}, \req{cgt54:temp-7} and \req{plarge:fct-d2} continue to hold in this case (they only use general ARp properties and the mechanism of
the algorithm). Applying \req{psmall:sigma-upper} in \req{plarge:fct-d2}, we deduce
\eqn{psmall:fct-d}{f(x_k)-f(x_{k+1})\geq \frac{\eta_1}{r}(\alpha\epsilon)^{\frac{r}{r-1}}\sigma_{\rm up}^{-\frac{1}{r-1}}= 
\frac{\eta_1}{r}\left(\frac{\alpha^r}{\sigma_{\rm up}}\right)^{\frac{1}{r-1}}\cdot\epsilon^{\frac{r}{r-1}},}
which gives \req{psmall:fct-final}. 

Using that $f(x_k)=f(x_{k+1})$ on unsuccessful iterations, and that $f(x_k)\geq f_{\rm low}$ for all $k$, we can sum up over all successful iterations to deduce 
\req{psmall:succ-no}.}



We are left with counting the number of  total iterations and evaluations.

\lcor{psmall:final-cor}{Let $r\leq p+\beta_p$, assume A.1 and that $\{f(x_k)\}$ is bounded below by $f_{\rm low}$. 
 If $r<p+\beta_p$ assume also A.2 and $\sigma_{\min}>0$. Then Algorithm \ref{ARq-overest} takes at most 
\eqn{psmall:total-no}{
\left\lfloor \frac{f(x_0) - f_{\rm low}}{\kappa_{s,r}}\left(1+ \frac{|\log\gamma_3|}{\log\gamma_1}\right) \epsilon^{-\frac{r}{r-1}} 
+ \frac{1}{\log\gamma_1}\log \frac{\sigma_{\rm up}}{\sigma_0}\right\rfloor
}
iterations/evaluations of $f$ and its derivatives until termination, where $\kappa_{s,r}$ and $\sigma_{\rm up}$ are defined in \req{psmall:succ-no} and \req{psmall:sigma-upper},
respectively.}

\proof{We first upper bound the total number of unsuccessful iterations; for this, we apply  Lemma \ref{it-count-total-0} to upper bound
$|\mathcal{U}_j|$ with $\sigma_{\rm up}$ defined in \req{psmall:sigma-upper}. To prove \req{psmall:total-no}, 
use \req{psmall:succ-no} and \req{Uj}, where we let $j$ denote the first iteration with $\pi_f(x_j+s_j)<\epsilon$
(so the iteration where ARp terminates), and we use $j=|\mathcal{S}_j|+|\mathcal{U}_j|$.}

\begin{remark}\label{remark4}
{\rm \begin{itemize}
\item[(a)] (Comment on $\sigma_{\min}$) Note that $\sigma_{\min}>0$ only appears/is used in the complexity bounds for the regime $r<p+\beta_p$ (namely in the
definition of the constant $M$ in A.2) and not for the case $r=p+\beta_p$ (see also our Remark \ref{remark3} (a)). 

\item[(b)] (Condition \req{step-long-gn}) We have used \req{step-long-gn} in the proof of Theorem 
\ref{p-small:thm-success}
(namely, in the use of \req{plarge:fct-d2} to deduce
\req{psmall:fct-d})
and hence for obtaining
the main complexity result in the regime $p<r\leq p+\beta_p$. This was however, not strictly necessary for obtaining same order
complexity bounds (albeit with different constants) in this parameter
regime, and was done for simplicity and coherence of the algorithm and results with the regime $r>p+\beta_p$ (for which 
\req{step-long-gn} is needed), and for practicality as $\beta_p$ is not known {\it a priori}.
 Let us briefly outline how one could bypass the use of \req{step-long-gn} in the proof of 
Theorem \ref{p-small:thm-success}. Note first that \req{step-long-gn} implies in this regime, given the constant upper bound 
\req{psmall:sigma-upper}, that $\|s_k\|\geq {\rm constant}\times \epsilon^{\frac{1}{r-1}}$. A similar lower bound on $s_k$ can be obtained 
directly (rather than from \req{step-long-gn}) from \req{step-lowerbd} as follows: when $\|s_k\|\leq 1$, 
 \req{step-lowerbd} implies $(\sigma_k+\theta+\kappa_2)\|s_k\|^{r-1}\geq \epsilon$; thus, using the constant upper bound \req{psmall:sigma-upper} on $\sigma_k$, $\|s_k\|\geq\min\{1, {\rm constant}_{\rm new}\times \epsilon^{\frac{1}{r-1}}\}$. Using the latter bound
 in \req{plarge:fct-d}, and that $\sigma_k\geq \sigma_{\min}$ and $\epsilon\in (0,1]$, we can deduce a same-order bound (in 
 $\epsilon$) as in \req{psmall:fct-final}. This line of proof is remindful of techniques used in \cite{bgmst} (for the case
 $\beta_p=1$ and $r=p+1$).
\item[(c)] (The Lipschitz continuous case) Letting $\beta_p=1$ (namely, the $p$th order derivative is Lipschitz continuous)
 and $r=p+1$  recovers the complexity bounds in 
\cite{bgmst}, namely, $\mathcal{O}\left(\epsilon^{-\frac{p+1}{p}}\right)$ (albeit with different constants), and
shows these bounds continue to hold for any $r\geq p+1$. Note however, that condition \req{step-long-gn} is not needed 
in the ARp algorithm in \cite{bgmst}. Our previous remark (b) explains that 
 \req{step-long-gn} is not strictly needed for the complexity bounds 
in the regime $r\leq p+\beta_p$ (which includes the case $\beta_p=1$ and $r=p+1$) for our ARp variant, 
which clarifies the connection with the algorithm in \cite{bgmst}.
\item[(d)] (The case $r=p+\beta_p$) Despite their different proofs, when $r=p+\beta_p$, 
 the complexity bound \req{psmall:total-no} is {\it identical} to the (limit of the) bound \req{total-no}. Comparing the expressions of these two bounds, we find that $r=p+\beta_p$ implies that the $|\log \epsilon|$ term  in  \req{total-no} vanishes, and that
the two complexity bounds clearly agree provided   $\kappa_{s,p}=\kappa_{s,r}$ and $\sigma_{\max}=\sigma_{{\rm up}}$. 
Furthermore, the definitions \req{constants-success} and \req{psmall:constants-success} trivially imply $\kappa_{s,p}=\kappa_{s,r}$ if $\sigma_{\max}=\sigma_{{\rm up}}$.
Finally, to see the latter identity, use the corresponding definitions in \req{constants-success} and \req{psmall:sigma-upper} and note that 
 $r=p+\beta_p$ provides that $\kappa_1=\kappa_2$, where $\kappa_1$ is defined in \req{constants}.
\end{itemize}}
\end{remark}

The constants in the bound \req{psmall:total-no} and their behaviour with respect to increasing values of $p$  are discussed in Section \ref{constants-page}.




\subsection{The constants in the complexity bounds}
\label{constants-page}

In this section we extract the key constants and expressions in the complexity bounds \req{total-no} and
\req{psmall:total-no} with respect to  $p$ and $r$ and show that in important cases, they stay finite as $p$ grows, for some suitable choices of algorithm parameters.

{\bf The case $r=p+1$, $\beta_p\in [0,1]$, $p\geq 2$.} \quad In this case, the complexity bound \req{total-no} applies
for $\beta_p\in [0,1)$. When $\beta_p=1$ (the Lipschitz continuous case), the bound \req{psmall:total-no} holds;
however, in Remark \ref{remark4} (d), we showed that 
\req{psmall:total-no} and (the limit of) \req{total-no} coincide when $r=p+\beta_p=p+1$. Hence, without loss of generality,
we focus on estimating \req{total-no} for any $\beta_p\in [0,1]$. Again without prejudice, we ignore
algorithm parameters (namely, $\gamma_1$, $\gamma_2$ and $\gamma_3$) that are independent of $p$ as they can easily be fixed.
Then, 
 \req{total-no} is a constant
 multiple of 
 \eqn{total-no-2}{
\left\lfloor \frac{f(x_0) - f_{\rm low}}{\kappa_{s,p}} \epsilon^{-\frac{p+\beta_p}{p+\beta_p-1}} 
+ \frac{(1-\beta_p)|\log \epsilon|}{p+\beta_p-1} + \log \frac{\sigma_{\rm max}}{\sigma_0}\right\rfloor.
}
From \req{constant-0} and \req{constants}, we deduce 
 \eqn{constants-approx-2}{
 \kappa_2=\mathcal{O}(L_p)\quad{\rm and}\quad \kappa_1 = 3^{\frac{1-\beta_p}{p+\beta_p-1}}\kappa_2^{\frac{p}{p+\beta_p-1}}=\mathcal{O}\left(L_p^{\frac{p}{p+\beta_p-1}}\right),
 }
 and hence, from \req{constants-success},
 \eqn{constants-approx}{
 \sigma_{\max}=\max\{\sigma_0,\gamma_2\theta, \gamma_2\kappa_1\}\quad {\rm and}\quad \frac{1}{\kappa_{s,p}}=\mathcal{O}\left((p+1)\sigma_{\max}^{\frac{1}{p}}\right)=  \mathcal{O}\left((p+1)\max\{\sigma_0^{\frac{1}{p}},
 \theta^{\frac{1}{p}}, L_p^{\frac{1}{p+\beta_p-1}}\}\right)
 }
 where we note that the term $(p+1)$ arises from the denominator of \req{2:localmodels} and $r=p+1$. 
Note that for simplicity of calculations, the H\"older constant $L_p$ in A.1 was scaled by $(p-1)!$. Thus letting
$L$ denote the usual/unscaled H\"older constant, we have 
\eqn{LipsConst}{
L\eqdef (p-1)!L_p,} 
where we assume that $L$ is independent, or stays
bounded with $p$. (Of course, $L$ and $L_p$ can have further
implicit dependencies on $p$ which are difficult to make precise.) 

Taking \req{LipsConst}
explicitly into account, and using 
Stirling's formula $\left[(p-1)!\sim [(p-1)/e]^{p-1}\sqrt{2\pi (p-1)}\right]$, we deduce 
\eqn{stirling}{
\begin{array}{lcl}\displaystyle
\lim_{p\rightarrow \infty} (p+1)L_p^{\frac{1}{p+\beta_p-1}}&=&\lim_{p\rightarrow \infty} (p+1)\left(\frac{L}{(p-1)!}\right)^{\frac{1}{p+\beta_p-1}}\\
&=&\displaystyle
\lim_{p\rightarrow \infty} (p+1)L^{\frac{1}{p+\beta_p-1}}[2\pi(p-1)]^{-\frac{1}{2(p+\beta_p-1)}}\left(\frac{p-1}{e}\right)^{-\frac{p-1}{p+\beta_p-1}}\\
&=&\displaystyle \lim_{p\rightarrow \infty} \left(\frac{L}{\sqrt{2\pi}}\right)^{\frac{1}{p+\beta_p-1}}\times
\lim_{p\rightarrow \infty} (p+1)(p-1)^{-\frac{1}{2(p+\beta_p-1)}}\left(\frac{p-1}{e}\right)^{-\frac{p-1}{p+\beta_p-1}}\\
&=&1\times\displaystyle
\lim_{p\rightarrow \infty}(p-1)^{-\frac{1}{2(p+\beta_p-1)}}e^{\frac{p-1}{p+\beta_p-1}}\frac{p+1}{(p-1)^{\frac{p-1}{p+\beta_p-1}}}  =1\times e\times 1 =e,
\end{array}
}
where we used the standard limits $\lim_{u\rightarrow\infty} u^{\frac{1}{u}}=1$ and $\lim_{u\rightarrow\infty} c^{\frac{1}{u}}=1$,
where $c>0$ is an arbitrary constant.
This and \req{constants-approx} imply that 
\[
\lim_{p\rightarrow \infty} \frac{1}{\kappa_{s,p}}<\infty,
\]
provided that 
\eqn{finite-lim}{
(p+1)\sigma_0^{\frac{1}{p}}<\infty \quad {\rm and}\quad (p+1)\theta^{\frac{1}{p}}<\infty, \quad{\rm as}\quad p\rightarrow \infty. 
}
The limits in \req{finite-lim} can be 
achieved without difficulty by suitable choices/scalings of  $\sigma_0$ and $\theta$, which  are user-chosen algorithm parameters. 
In particular,   let 
\eqn{sigma0-theta}{
\sigma_0\eqdef \frac{\overline{\sigma}_0}{(p-1)!} \quad {\rm and}\quad \theta\eqdef \frac{\overline{\theta}}{(p-1)!}, 
}
for any constants $\overline{\sigma}_0$ and $\overline{\theta}$ independent of $p$; Stirling's formula applied to $(p-1)!$
and similar calculations to \req{stirling} can be used to show that \req{sigma0-theta} satisfy \req{finite-lim}.

The second term in the sum \req{total-no-2} either vanishes when $\beta_p=1$ or converges to zero as $p\rightarrow 0$.
Proceeding to the third term in the sum \req{total-no-2}, we have:
from  \req{constants-approx-2} and \req{LipsConst}, we deduce $\kappa_1\rightarrow 0$ as $p\rightarrow \infty$ and so, irrespective 
of the scaling of $\sigma_0$ and $\theta$,
$1\leq \sigma_{\max}/\sigma_0<\infty$. Thus the last term in \req{total-no-2} is finite.
 
We can safely conclude now that as $p\rightarrow \infty$, all constants in \req{total-no-2}
stay bounded or converge to zero for appropriate choices of $\sigma_0$ and $\theta$, and so, using also that $\epsilon \in (0,1]$, the bound \req{total-no} approaches $\mathcal{O}(\epsilon^{-1})$.

The above discussion of limiting constants can be easily extended,
with similar results, to any $r=ap+b$
with $a, b>0$ independent of $p$, provided $r>p+\beta_p$.

Note also that the more practical case is when $p$ is fixed and $\epsilon$ can be made arbitrarily small; then, the bound 
\req{total-no} is well-defined for all algorithm and problem parameter choices, allowing the use of simplified constants and unscaled
parameters in the analysis.\\[1ex]


 
 
 
{\bf The case $r=p+\beta_p$, $\beta_p\in [0,1]$, $p\geq 2$.}\quad In this case, the bound \req{psmall:total-no} applies
(note that the case $\beta_p=1$ was already addressed in the first case of this section). 
The constants in \req{psmall:total-no}  stay
bounded as $p$ grows,  provided  $\sigma_0$ and $\theta$ are scaled according to \req{sigma0-theta}. Indeed, one can show this
very similarly to the case $r=p+1$ above, using \req{constant-0}, \req{psmall:constants-success} and \req{LipsConst} to obtain the following estimates
\[
\kappa_2=\mathcal{O}(L_p)=\mathcal{O}\left(\frac{L}{(p-1)!}\right),\quad \sigma_{\rm up}=\max\{\sigma_0,\gamma_2\theta,\gamma_2\kappa_2\}=\mathcal{O}(\max\{\sigma_0,\theta, L_p\}).
\]
Letting $r=p+\beta_p$ in \req{psmall:constants-success}, we have
\[
\frac{1}{\kappa_{s,r}}=\mathcal{O}\left(r\sigma_{\rm up}^{\frac{1}{r-1}}\right)=\mathcal{O}\left((p+\beta_p)\sigma_{\rm up}^{\frac{1}{p+\beta_p-1}}\right)=\mathcal{O}\left((p+\beta_p)(\max\{\sigma_0,\theta, L_p\})^{\frac{1}{p+\beta_p-1}}\right)<\infty, \quad{\rm as}\quad p\rightarrow\infty,
\]
where the limit follows similarly to \req{stirling}, using also \req{sigma0-theta}.
As $p$ grows and as
a function of $\epsilon$, \req{psmall:total-no} approaches
the same well-defined limit as \req{total-no}, namely, $\mathcal{O}(\epsilon^{-1})$.\\[1ex]


{\bf The case $p<r<p+\beta_p$, $\beta_p\in [0,1]$, $p\geq 2$.}\quad In this case, the bound \req{psmall:total-no} applies. 
  However, the limiting constants
in \req{psmall:total-no}
 depend crucially on $M$ in A.2, which grows unbounded with $p$.


\numsection{Discussion of complexity bounds}

\subsection{The cubic regularization algorithm}
\label{Univ-condition-explanation}

We now particularize our algorithm and results to the case when $p=2$ and $r=p+1$,
which yields a cubic regularization model \req{2:localmodels} and algorithm, with condition \req{step-long-gn}, namely,
\eqn{cubic:step-long-gn}{
\sigma_k\|s_k\|^2\geq \alpha \pi_f(x_k+s_k),}
imposed on any successful step $s_k$, and which allows $\sigma_{\min}=0$ in \req{sigupdate}.

\lcor{cubic:final}{Let $p=2$, $r=3$ and $\epsilon \in (0,1]$.
 Assume that $f\in C^2(\mathcal{F})$, and $\nabla^2_x f$ is H\"older continuous on the path of the iterates
 and trial points with exponent $\beta_2\in [0,1]$. Let $\{f(x_k)\}$ be bounded below by $f_{\rm low}$. 
Then for all successful iterations $k$ until the termination of
Algorithm  \ref{ARq-overest}, we have
\eqn{cubic:fct-final}{f(x_k)-f(x_{k+1})\geq \kappa_{s,2} \epsilon^{\frac{2+\beta_2}{1+\beta_2}},}
where 
\eqn{cubic:constants-success}{
\kappa_{s,2}\eqdef \frac{\eta_1}{3}\left(\frac{\alpha^{3}}{\sigma_{\max}}\right)^{\frac{1}{2}},\quad \sigma_{\max}\eqdef \max\left\{\sigma_0,\gamma_2\theta, \gamma_2\kappa_1\right\},
}
and $\kappa_1\eqdef 3^{\frac{3-\beta_2}{1+\beta_2}}\left[\frac{L_2}{2(1-\eta_2)}\right]^{\frac{2}{1+\beta_2}}$. Thus Algorithm \ref{ARq-overest} takes at most 
\eqn{cubic:succ-no}{
\left\lfloor \frac{f(x_0) - f_{\rm low}}{\kappa_{s,2}} \epsilon^{-\frac{2+\beta_2}{1+\beta_2}}\right\rfloor
}
successful iterations/evaluations of derivatives of degree $2$  of $f$ until termination,
and 
at most 
\eqn{cubic:total-no}{
\left\lfloor \frac{f(x_0) - f_{\rm low}}{\kappa_{s,2}}\left(1+ \frac{|\log\gamma_3|}{\log\gamma_1}\right) \epsilon^{-\frac{2+\beta_2}{1+\beta_2}} 
+ \frac{1-\beta_2}{(1+\beta_2)\log\gamma_1}|\log \epsilon| + \frac{1}{\log\gamma_1}\log \frac{\sigma_{\rm max}}{\sigma_0}\right\rfloor
}
iterations/evaluations of $f$ and its first and second derivatives until termination, where $\kappa_{s,2}$ and $\sigma_{\max}$ are defined in \req{cubic:constants-success}.}

\proof{Clearly, the results follow from Corollary \ref{cor-complexity} for $p=2$, $r=3$ and $\beta_2\in [0,1)$, and from Corollary \ref{psmall:final-cor}
for $p=2$, $r=3$ and $\beta_2=1$. We note the key ingredients that are needed to obtain \req{cubic:fct-final},
with the remaining results following from standard telescopic sum arguments and from Lemma \ref{it-count-total-0}, respectively. Lemmas 
\ref{lem-sigma-upper-bd} and \ref{p:small-s-upper-bd} provide the following upper bound on $\sigma_k$,
\[
\sigma_k\leq \sigma_{\max}\epsilon^{-\frac{1-\beta_2}{1+\beta_2}}, \quad k\geq 0.
\]
This bound and  condition \req{cubic:step-long-gn} (which is \req{step-long-gn}) are then substituted into 
the objective decrease condition \req{plarge:fct-d} on successful steps which here takes the form
\[
f(x_k)-f(x_{k+1})\geq \frac{\eta_1}{3}\sigma_k\|s_k\|^3\geq \frac{\eta_1}{3}\alpha\epsilon\left(\frac{\alpha\epsilon}{\sigma_k}\right)^{\frac{1}{2}}\geq \frac{\eta_1}{3}\left(\frac{\alpha^3}{\sigma_{\max}}\right)^{\frac{1}{2}}\epsilon^{\frac{3}{2}}.
\]
}

The impact of the value of $\beta_2\in [0,1]$ can be seen in the bound \req{cubic:total-no}; for example, when $\beta_2=1$, the
$|\log\epsilon|$ term disappears, in agreement with known bounds for ARC \cite{ARCpartII}. Note that as a function of $\epsilon$,
Corollary \ref{cubic:final} matches corresponding bounds in \cite{GN16} (for different cubic regularization variants)  and extends them to convex constraints, allowing inexact subproblem solves. Our purpose here is also to allow $p\geq 2$, and a discussion of the bounds we obtained follows.



\subsection{General discussion of the complexity bounds}

Table \ref{table1} gives a summary of our complexity bounds as a function of $r$ and $q$. 

\begin{table}\label{table1}
\begin{center}
\begin{tabular}{|l|c|c|}
\hline
      Algorithm                         &        $p<r\leq p+\beta_p $  & $p+\beta_p<r$ \\[0.25ex]
\hline
ARp with $p=1$  & $\mathcal{O}\left(\epsilon^{-\frac{r}{r-1}}\right)=\left[\mathcal{O}\left(\epsilon^{-\frac{1+\beta_1}{\beta_1}}\right),\infty\right)$ &   $\mathcal{O}\left(\epsilon^{- \frac{1+\beta_1}{\beta_1}}\right)$\\*[2ex]
\hline
ARp with $p=2$     &  $\mathcal{O}\left(\epsilon^{-\frac{r}{r-1}}\right)=\left[\mathcal{O}\left(\epsilon^{-\frac{2+\beta_2}{1+\beta_2}}\right),\mathcal{O}\left(\epsilon^{-2}\right)\right)$  & $\mathcal{O}\left(\epsilon^{-\frac{2+\beta_2}{1+\beta_2}}\right)$  \\*[2ex]
\hline
ARp with $p=3$     &   $\mathcal{O}\left(\epsilon^{-\frac{r}{r-1}}\right)=\left[\mathcal{O}\left(\epsilon^{-\frac{3+\beta_3}{2+\beta_3}}\right),\mathcal{O}\left(\epsilon^{-\frac{3}{2}}\right)\right)$  & $\mathcal{O}\left(\epsilon^{-\frac{3+\beta_3}{2+\beta_3}}\right)$  \\*[2ex]
\hline
$\ldots$&$\ldots$&$\ldots$\\*[2ex]
\hline
ARp with $p\geq 2$  &  $\mathcal{O}\left(\epsilon^{-\frac{r}{r-1}}\right)=\left[\mathcal{O}\left(\epsilon^{-\frac{p+\beta_p}{p+\beta_p-1}}\right),\mathcal{O}\left(\epsilon^{-\frac{p}{p-1}}\right)\right)$  & $\mathcal{O}\left(\epsilon^{-\frac{p+\beta_p}{p+\beta_p-1}}\right)$  \\*[2ex]
\hline
\end{tabular}
\end{center}
\caption{Summary of complexity bounds for regularization methods for ranges of $r$. Recall we assumed that $\epsilon \in (0,1]$, $r>p\geq 1$, $r\in \Re$ and $p\in \Na$; and that either $p\geq 1$ and $\beta_p\in (0,1]$, or $p\geq 2$ and $\beta_p\in [0,1]$. Also, the ranges in the second column are as a function of the dominating terms in $\epsilon$ and varying $r$ in the appropriate interval and they are 
plotting the changing bound $\mathcal{O}(\epsilon^{\frac{r}{r-1}})$.}
\end{table}


Several remarks and comparisons are in order concerning these bounds.
\begin{itemize}
\item {\bf The first-order case.} Note that the case $p=1$ is also covered,  with a more general quadratic model and 
using a Cauchy analysis, in \cite{cgt51}; the same complexity bounds ensue (as a function of the
accuracy) as in Table \ref{table1} for $p=1$; the case $\beta_1=0$ is also not covered in \cite{cgt51}. 

\item {\bf Sharpness.} For unconstrained problems ($\mathcal{F}=\Re^n$), the bound for the case $p=1$ and $r\geq 1+\beta_1$, $\beta_1\in (0,1]$, 
was shown to be sharp in \cite{cgt51}. 
Also, the bounds for ARp with $p=2$ and $2<r\leq 2+\beta_2$ and $\beta_2\in (0,1]$ are sharp and optimal
for the corresponding smoothness classes \cite{cgt38}. We also note that for general $p$, $r=p+1$ and
$\beta_p=1$ (the Lipschitz continuous case),  \cite{Carmon} shows the bounds  for (possibly randomized)
ARp variants (in \cite{bgmst}) 
to be sharp and optimal. The difficult example functions in \cite{Carmon} increase in dimension with $p$,
in contrast to uni- or bi-variate examples in \cite{cgt51, cgt38}.

\item {\bf Continuity.} All bounds vary continuously with $r$ and $\beta_p\in [0,1]$. In particular, when $r= p+\beta_p$, the complexity bounds in the second and third column match
(for a given $p$ and $\beta_p$) (see also Remark \ref{remark4} (d)).  

\item {\bf Universality \cite{NemirovskiYudin, Nest13b, GN16}.} 
For fixed $p$ and $\beta_p$, the best complexity bounds are obtained when $r\geq p+\beta_p$. 
These bounds do not depend on the regularization power $r$, and even though the smoothness parameter $\beta_p$ is (usually) unknown, its value is captured accurately in the
complexity, even for the case when $\beta_p=0$ and  $p\geq 2$.
Note that the values of the complexity bounds as a function of the accuracy indicate that one
should choose $r\geq p+1$ to achieve the best complexity when
$\beta_p$ is unknown; and there
seems to be little reason, from an evaluation complexity point of
view, to pick anything other than $r= p+1$. (But, note that, as a benefit of using \req{step-long-gn}, one can simplify ARp's construction by not imposing 
a lower bound $\sigma_{\min}$ in the $\sigma_k$ update \req{sigupdate}.)

\item {\bf Complexity values in the order of the accuracy.}
Table \ref{table1} shows the increasingly good complexity obtained as $p$ grows and $\beta_p\in [0,1]$,
namely, the more derivatives are available and the smoother these derivatives are. In particular, purely as a function of $\epsilon$ and as $r$ varies, we obtain the following ranges of complexity
powers : $[\epsilon^{-2},\infty)$ ($p=1$); $[\epsilon^{-\frac{3}{2}},\epsilon^{-2}]$ ($p=2$); $[\epsilon^{-\frac{4}{3}},\epsilon^{-\frac{3}{2}}]$ ($p=3$); $[\epsilon^{-\frac{5}{4}},\epsilon^{-\frac{4}{3}}]$ ($p=4$);
and so on. 

\item {\bf The Lipschitz continuous case.} Letting $\beta_p=1$ (namely, the $p$th order derivative is Lipschitz continuous)
 and $r=p+1$ in Table \ref{table1} recovers the complexity bounds in 
\cite{bgmst}, namely, $\mathcal{O}\left(\epsilon^{-\frac{p+1}{p}}\right)$; see also  Remark \ref{remark4} (c).
Furthermore, the results here show that for our ARp variant, this complexity bound continues to hold
for any regularization power $r\geq p+1$.

\item {\bf Loss of smoothness}  Note that for fixed $p\geq 2$, $\beta_p=0$  corresponds to the case when the objective has the highest level of non-smoothness compared to $\beta_p\in (0,1]$. 
Then ARp can still be applied, and the good complexity bounds for the case $r\geq p+\beta_p\geq 2$ hold.

\item {\bf Constants in the complexity bounds} The constants in the complexity bounds for $r\geq p+\beta_p$ stay bounded (above)
as $p$ grows, provided some user-chosen algorithm parameters are suitably scaled and that $r=O(p)$ (see Section \ref{constants-page}). Thus these complexity
bounds remain valid with growing $p$ and approach $\mathcal{O}(\epsilon^{-1})$.

\end{itemize}

\numsection{Conclusions}

We have generalized and modified the regularization methods in \cite{bgmst} to allow for varying regularization power, accuracy of Taylor polynomials and different (H\"older) smoothness levels
of derivatives. Our results show the robustness of the evaluation complexity bounds with respect to such perturbations. We found that complexity bounds of regularization methods improve
with growing accuracy of the Taylor models and increasing smoothness levels of the objective. Furthermore, when the regularization power $r$ is sufficiently large
(say $r\geq p+1$) our modification to ARp in the spirit of \cite{GN16} 
allows ARp's worst-case behaviour to be independent of the regularization power and to accurately reflect the (often unknown) smoothness level of the objective.
We have also generalized \cite{bgmst} and \cite{GN16}  to problems with convex constraints and inexact subproblem solutions.
The question as to whether the complexity bounds we obtained are sharp remains open when $r\neq p+\beta_p$ and $p\geq 3$.
This question is particularly poignant in the case when $p<r<p+\beta_p$:
could a suitable modification of ARp achieve an (improved) evaluation complexity bound that is independent of the regularization power in this case as well?


\begin{thebibliography}{10}

\bibitem{BensFreh02}
Alain Bensoussan and Jens Frehse.
\newblock {\em Regularity results for nonlinear elliptic systems and
  applications}.
\newblock Springer Verlag, Heidelberg, Berlin, New York, 2002.

\bibitem{Bertsekas}
D.P.Bertsekas.
\newblock {\em Nonlinear Programming}. 
\newblock Athena Scientific, Belmont, Massachusetts, USA, 2nd edition, 1999.

\bibitem{bgmst}
E.~G. Birgin, J.~L. Gardenghi, J.~M. Mart\'{i}nez, S.~A. Santos, and Ph.~L.
  Toint.
\newblock Worst-case evaluation complexity for unconstrained nonlinear
  optimization using high-order regularized models.
\newblock {\em Mathematical Programming, Series A}, 163(1--2):359--368, 2017.

\bibitem{Birgin}
E.G. Birgin and  J.M. Mart\'{i}nez.
\newblock The use of quadratic regularization with a cubic descent condition for unconstrained optimization.
\newblock {\em SIAM Journal on Optimization}, 27(2):1049--1074, 2017.

\bibitem{bgms1}
E.~G. Birgin, J.~L. Gardenghi, J.~M. Mart\'{i}nez, and S.~A. Santos.
\newblock Remark on Algorithm 566: Modern Fortran routines for testing unconstrained optimization software with derivatives up to third-order.
\newblock Technical report, Department of Computer Science, University of S\~ao Paulo, Brazil, 2018.

\bibitem{bgms2}
E.~G. Birgin, J.~L. Gardenghi, J.~M. Mart\'{i}nez, and S.~A. Santos.
\newblock On the use of third-order models with fourth-order regularization for unconstrained optimization.
\newblock Technical report, Department of Computer Science, University of S\~ao Paulo, Brazil, 2018.

\bibitem{Carmon}
Y. Carmon, J.C. Duchi, O. Hinder and A. Sidford.
\newblock Lower bounds for finding stationary points I.
\newblock Technical report, arXiv:1710.11606, 2017. 

\bibitem{cgt36}
C.~Cartis, N.~I.~M. Gould, and Ph.~L. Toint.
\newblock On the complexity of steepest descent, {N}ewton's and regularized
  {N}ewton's methods for nonconvex unconstrained optimization.
\newblock {\em SIAM Journal on Optimization}, 20(6):2833--2852, 2010.


\bibitem{ARCpartII}
C.~Cartis, N.~I.~M. Gould, and Ph.~L. Toint.
\newblock Adaptive cubic overestimation methods for unconstrained optimization.
  {P}art {II}: worst-case function-evaluation complexity.
\newblock {\em Mathematical Programming, Series~A}, 130(2):295--319, 2011.

\bibitem{cgt38}
C.~Cartis, N.~I.~M.~Gould and Ph.~L.~Toint.
\newblock Optimal Newton-type methods for nonconvex smooth optimization problems.
\newblock ERGO Technical Report 11-009, School of Mathematics, University of Edinburgh, 2011.

\bibitem{cgt51}
C.~Cartis, N.~I.~M.~Gould and Ph.~L.~Toint.
\newblock Worst-case evaluation complexity of regularization methods for smooth unconstrained optimization using
{H}older continuous gradients.
\newblock {\em Optimization Methods and Software}, 32(6):1273--1298, 2017.

\bibitem{Chen}
X.~Chen, Ph.~L.~Toint and H.~Wang.
\newblock Partially separable convexly-constrained optimization with non-Lipschitzian singularities and its complexity.
\newblock Technical report, arXiv 1704.06919, 2017.

\bibitem{TRbook}
A.R. Conn, N.I.M. Gould and Ph. L. Toint.
\newblock {\em Trust region methods}.
\newblock MOS-SIAM series on Optimization, 2000.

\bibitem{TRACE}
F.~E. Curtis, D.~P. Robinson, and M.~Samadi.
\newblock A trust region algorithm with a worst-case iteration complexity of
  {$O(\epsilon^{-3/2})$} for nonconvex optimization.
\newblock {\em Mathematical Programming, Series~A}, 162(1):1--32, 2017.


\bibitem{Devo13}
O.~Devolder.
\newblock Exactness, Inexactness and Stochasticity in First-Order Methods for Large-Scale Convex Optimization.
\newblock PhD Thesis, ICTEAM and CORE, Universit\'e Catholique de Louvain, 2013.


\bibitem{GPSA94}
Gas Processors and Suppliers Association.
\newblock {\em Engineering Data Book. Vol.~2}.
\newblock GPSA, Tulsa, USA, 1994.

\bibitem{Gould1}
N.I.M. Gould, T. Rees and J. Scott.
\newblock A higher-order method for solving nonlinear least-squares problems.  
\newblock RAL Preprint RAL-P-2017-010, STFC Rutherford Appleton Laboratory, 2017.


\bibitem{Gould2}
N.I.M. Gould, T. Rees and J. Scott.
\newblock Convergence and evaluation-complexity analysis of a regularized tensor-Newton method for solving nonlinear least-squares problems. 
\newblock RAL Preprint RAL-P-2017-009, STFC Rutherford Appleton Laboratory, 2017.


\bibitem{GN16}
G.~N.~Grapiglia and Yu.~Nesterov.
\newblock Globally-convergent second-order schemes for minimizing twice-differentiable functions.
\newblock {\em SIAM Journal on Optimization}, 27(1):478--506, 2017.

\bibitem{Grie81}
A.~Griewank.
\newblock The modification of {N}ewton's method for unconstrained optimization
by bounding cubic terms.
\newblock Technical Report NA/12 (1981), Department of Applied Mathematics and Theoretical Physics,
University of Cambridge, United Kingdom, 1981. 

\bibitem{NemirovskiYudin}
A.~S.~Nemirovski and D.~B.~Yudin.
\newblock {\em Problem Complexity and Method Efficiency in Optimization}
\newblock Wiley Interscience Series in Discrete Mathematics, 1983.


\bibitem{Nest04}
{Yu}. Nesterov.
\newblock {\em Introductory Lectures on Convex Optimization}.
\newblock Applied Optimization. Kluwer Academic Publishers, Dordrecht, The
  Netherlands, 2004.


\bibitem{Nest13b}
{Yu}. Nesterov.
\newblock Universal gradient methods for convex optimization problems.
\newblock  {\em Mathematical Programming, Series~A}, 152(1--2):381--404, 2015.

\bibitem{NestModels}
{Yu}. Nesterov.
\newblock Implementable tensor methods in unconstrained convex optimization. 
\newblock CORE Discussion Paper, Universite Catholique de Louvain, Belgium, 2015.

\bibitem{NestPoly06}
{Yu}. Nesterov and B.~T. Polyak.
\newblock Cubic regularization of {N}ewton method and its global performance.
\newblock {\em Mathematical Programming, Series~A}, 108(1):177--205, 2006.

\bibitem{Schnabel}
R.B. Schnabel and T.T. Chow.
\newblock Tensor methods for unconstrained optimization using second derivatives.
\newblock {\em SIAM Journal on Optimization}, 1(3):293--315, 1991.

\bibitem{WeisDeufErdm07}
M.~Weiser, P.~Deuflhard and B.~Erdmann.
\newblock Affine conjugate adaptive {N}ewton methods for nonlinear elastomechanics.
\newblock {\em Optimization Methods and Software}, 22(3):413--431, 2007. 


\end{thebibliography}
\end{document}